\documentclass[a4paper]{article}
\usepackage{fullpage}

\usepackage[T1]{fontenc}
\usepackage[utf8]{inputenc}
\usepackage{lmodern}
\usepackage{booktabs}
\usepackage{longtable}
\usepackage{array}
\usepackage{pdflscape}
\usepackage{hyperref}
\usepackage{siunitx}
\hypersetup{
	colorlinks = true,
	linkcolor  = blue,
	citecolor  = blue,
	urlcolor   = blue,
	pdftitle   = {Dynamic PPA for unconstrained minimization}
}

\usepackage{amsmath,amssymb}
\usepackage{graphicx}
\usepackage{bm,xcolor}
\usepackage{dsfont}
\usepackage{comment}
\usepackage{algorithm}
\usepackage{algpseudocode}

\usepackage{acronym}
\acrodef{PPA}{Proximal Point Algorithm}

\newcommand{\X}{{\color{red!60!yellow!30!black}\bm{x}}}

\newcommand{\XB}{\overline{\X}}

\newcommand{\Z}{{\color{red!60!yellow!30!black}\bm{z}}}

\newcommand{\Dx}{{\color{blue!60!yellow}\bm{D}_x}}
\newcommand{\cx}{{\color{blue!60!yellow}c_x}}
\newcommand{\cxmin}{{\color{blue!60!yellow}c_{x,\min}}}
\newcommand{\cxmax}{{\color{blue!60!yellow}c_{x,\max}}}

\newcommand{\diag}[1]{\mathrm{diag}\left(#1\right)}

\newcommand{\emailLink}[1]{\href{mailto:#1}{#1}}
\newcommand{\orcidLink}[1]{\href{https://orcid.org/#1}{#1}}
\newcommand{\email}[1]{\textsc{email} \emailLink{#1}}
\newcommand{\orcid}[1]{\textsc{orcid} \orcidLink{#1}}

\newcommand{\TheTitle}{Dynamic Proximal Point Method for Unconstrained Minimization}
\newcommand{\TheFunding}{This study was carried out within th project \emph{``SPACE IT UP! ASI contract n. 2024-5-E.0 CUP Master n. I53D24000060005''} funded by the Italian Space Agency (ASI), and the Ministry of University and Research (MUR).}

\begin{document}
	
\title{\bfseries\TheTitle}
\author{Enrico Bertolazzi\thanks{University of Trento, Via Sommarive 9, Trento, Italy, \email{enrico.bertolazzi@unitn.it}, \orcid{0000-0003-0487-5210}.}\and%
	Alberto De Marchi\thanks{University of the Bundeswehr Munich, Germany, \email{alberto.demarchi@unibw.de}, \orcid{0000-0002-3545-6898}.}\and%
	Davide Stocco\thanks{University of Trento, Via Sommarive 9, Trento, Italy, \email{davide.stocco@unitn.it}, \orcid{0000-0001-5902-4048}.}}
\date{}

\maketitle

\begin{abstract}
	In this work, we present a novel dynamic proximal point algorithm for unconstrained optimization. The method generates a sequence of proximal subproblems, where the quadratic regularization term is weighted by a diagonal matrix that is updated adaptively at each iteration. Each subproblem is solved using an inner Newton's method combined with a line search, which provides a global convergence mechanism for the nonlinear solver. At the outer level, the algorithm updates the reference point and adjusts the regularization parameter based on the performance of the inner Newton solver. We derive the reduced linear system used to compute the Newton step, define the corresponding merit function, and discuss practical approaches for constructing the diagonal scaling matrix from derivative information. The paper also provides implementation-oriented pseudocode and stopping criteria that are consistent with the proposed method.
\end{abstract}

\section{Introduction}

Unconstrained optimization is a fundamental problem in numerical analysis, scientific computing and engineering. It arises in a wide range of applications, including machine-learning model training, parameter estimation, inverse problems and the numerical solution of nonlinear models. Given a sufficiently smooth objective function, the goal is to compute a stationary point, preferably corresponding to a local, or in special cases global, minimizer. Newton-type methods are among the most effective algorithms for smooth unconstrained minimization. Near a non-degenerate minimizer, they typically exhibit fast local convergence. However, their global behavior can be unreliable when the Hessian is singular, ill-conditioned, indefinite or strongly affected by non-convexity. In such cases the Newton system may be poorly posed, and the computed Newton direction may fail to be a descent direction, possibly leading to instability or divergence~\cite{Nocedal2006}. Classical globalization strategies address these issues by modifying the local quadratic model or the Newton step. Trust-region methods restrict the step to a neighborhood where the model is expected to be reliable, whereas Levenberg--Marquardt-type methods regularize the Hessian by adding a positive multiple of the identity matrix~\cite{Levenberg1944,Marquardt1963,More1978}. Both of these approaches improve the conditioning of the linear system and promote descent. Another stabilization mechanism is provided by the \ac{PPA}, originally developed in convex analysis and nonsmooth optimization~\cite{Rockafellar1976,Parikh2014}. The proximal framework replaces the original problem by a sequence of regularized subproblems, each of which penalizes deviations from a prescribed anchor point, typically identified with the previous iterate. Such a regularization improves stability and controls the displacement between successive iterates.

In this work, we propose a \textit{Dynamic Proximal Point Method} for smooth unconstrained optimization. The method combines proximal regularization with curvature-aware diagonal scaling. Instead of using only an isotropic scalar regularization, the proximal term contains an adaptive diagonal weighting matrix constructed from local derivative information. This allows the regularization to reflect the anisotropic structure of the objective function. The resulting algorithm has a natural two-level structure: An \emph{outer} iteration updates the anchor point and adapts the regularization parameters, while an \emph{inner} Newton iteration solves the corresponding proximal subproblem. The method dynamically balances between gradient-like behavior, useful far from a minimizer or in ill-conditioned regions, and Newton-like behavior, effective near well-conditioned stationary points. This adaptivity improves robustness on non-convex and ill-conditioned problems while retaining the efficiency of second-order information when it is reliable.

\section{Problem Formulation}

We consider the unconstrained minimization problem
\begin{equation}
	\mathop{\textrm{minimize}}_{\X\in\mathbb{R}^n} \quad f(\X) ,
	\label{eq:original-problem}
\end{equation}
where $f:\mathbb{R}^n\to\mathbb{R}$ is twice continuously differentiable. The proposed method replaces~\eqref{eq:original-problem} by a sequence of proximal subproblems. Given an anchor point $\XB_\ell$, a positive diagonal scaling matrix $\Dx$ and a scalar regularization parameter $\cx>0$, we define the proximal merit function as
\begin{equation}
	\mathcal{M}_{\ell}(\X) = f(\X) + \dfrac{\cx}{2}\|\X-\XB_\ell\|_{\Dx}^2 .
	\label{eq:proximal-merit}
\end{equation}
The weighted norm $\|\cdot\|_{\Dx}$ is defined as
\begin{equation}
	\|\X\|_{\Dx}^2 = \X\cdot(\Dx\X) \quad \text{where} \quad \Dx = \diag{\bm{w}} \quad \text{and} \quad \bm{w}>0 .
	\label{eq:weighted-norm}
\end{equation}
Here, $\bm{w}$ provides component-wise scaling, while $\cx$ controls the overall strength of the proximal regularization.
Both quantities are updated dynamically during the outer iteration.

At outer iteration $\ell$, the next iterate is obtained by approximately solving
\begin{equation}
	\X_{\ell+1} \approx \mathop{\arg\min}_{\X\in\mathbb{R}^n} \mathcal{M}_{\ell}(\X) .
	\label{eq:outer-subproblem}
\end{equation}
The proximal term stabilizes the subproblem by improving the conditioning of the Hessian of the merit function and, when the regularization is sufficiently strong, by enforcing local convexity of the regularized model.

The algorithm can therefore be organized into two nested levels.
The \textbf{outer} loop updates the anchor point $\XB_\ell$, the diagonal scaling matrix $\Dx$, and the scalar regularization parameter $\cx$. The \textbf{inner} loop approximately minimizes the current proximal merit function $\mathcal{M}_{\ell}$. A high-level description is given in Algorithm~\ref{alg:algo-base}.

\begin{algorithm}[!htb]
	\caption{Dynamic proximal-point minimization algorithm}
	\label{alg:algo-base}
	\begin{algorithmic}[1]
		\State Choose an initial point $\X_0$
		\State Set the initial anchor point $\XB_0 \gets \X_0$
		\State Initialize the diagonal scaling matrix $\Dx$ and the scalar weight $\cx$ (Section~\ref{sec:Dx})
		\For{$\ell = 0,1,2,\ldots,\textrm{\{maximum number of outer iterations\}}$}
		\State Use Algorithm~\ref{alg:algo-inner} to solve the proximal subproblem (Section~\ref{sec:newton})
		\begin{equation*}
			\X_{\ell+1}
			\approx
			\mathop{\arg\min}_{\X\in\mathbb{R}^n}\;\mathcal{M}_{\ell}(\X)
		\end{equation*}
		\If{\textrm{\{the inner minimization fails\}}}
		\State Reject the trial point and set $\X_{\ell+1}\gets \X_{\ell}$
		\State Update the regularization weight $\cx$ (Section~\ref{sec:cx:reject})
		\Else
		\State Update the anchor point $\XB_{\ell+1}\gets \X_{\ell+1}$
		\State Update the diagonal scaling matrix $\Dx$ (Section~\ref{sec:Dx})
		\State Update the regularization weight $\cx$ (Section~\ref{sec:cx})
		\EndIf
		\State Check convergence (Section~\ref{sec:check:convergence}) and return if detect \textbf{converged} or \textbf{stagnate}.
		\EndFor
		\State Failure, reached maximum iteration limits.
	\end{algorithmic}
\end{algorithm}

\subsubsection{Adaptive Diagonal Scaling and Convergence Properties}
\label{sec:Dx}

The diagonal scaling matrix is defined by
\begin{equation}
	\Dx = \diag{\bm{w}}, \quad \text{with} \quad w_i = \max\left\{ R_i-D_i,\, \varepsilon_{\min} \right\} \quad \text{for} \quad i = 1,\ldots,n,
	\label{eq:adaptive-diagonal-scaling}
\end{equation}
and where
\begin{equation}
	D_i = \left(\nabla^2 f(\X)\right)_{ii} \quad \text{and} \quad R_i = \sum_{j=1}^{n} \left| \left(\nabla^2 f(\X)\right)_{ij} \right|.
	\label{eq:row-sum-scaling}
\end{equation}
Thus $R_i - D_i$ represents, up to the diagonal sign convention, the absolute off-diagonal row sum of the Hessian, while $\varepsilon_{\min} > 0$ prevents vanishing weights and keeps $\Dx$ positive
definite.

Let
\begin{equation*}
	\bm{H}(\X)=\nabla^2 f(\X) \quad \text{with} \quad \bm{J}(\X)=\bm{H}(\X)+\cx\Dx
\end{equation*}
be the regularized Hessian used in the inner Newton system. Since the
regularization is diagonal,
\begin{equation*}
	J_{ii}=H_{ii}+\cx w_i ,
	\quad \text{and} \quad
	J_{ij}=H_{ij}
	\quad \text{for} \quad
	\quad i \ne j .
\end{equation*}
The role of $\Dx$ is to increase each diagonal entry of $\bm{J}(\X)$ in proportion to the off-diagonal coupling in the corresponding row of
$\bm{H}(\X)$, thereby promoting strict diagonal dominance. Indeed, by the Gershgorin circle theorem, every eigenvalue of $\bm{J}(\X)$ belongs to some disk
\begin{equation*}
	\mathcal{G}_i = \left\{ z\in\mathbb{C} : |z-J_{ii}| \le \sum\nolimits_{j\ne i}|H_{ij}| \right\} .
\end{equation*}
Hence, if $\cx$ is large enough that
\begin{equation}
	H_{ii}+\cx w_i > \sum\nolimits_{j\ne i}|H_{ij}| \quad \text{for} \quad i = 1,\ldots,n ,
	\label{eq:strict-diagonal-dominance}
\end{equation}
then all Gershgorin disks lie in the positive real half-line. For symmetric $\bm{H}(\X)$, the regularized Hessian $\bm{J}(\X)$ is therefore symmetric positive definite, so the Newton direction is well defined and is a descent direction for the proximal merit function. This curvature-aware scaling requires no eigenvalue computations. If $\Dx$ is uniformly positive definite and $\mathcal{M}(\X)$ is twice continuously differentiable, increasing $\cx$ after failed inner solves and decreasing it after efficient solves gives a globalization mechanism analogous to trust-region or Levenberg--Marquardt regularization.

\section{Inner and Outer Iterations}

The proposed method has a nested structure. For fixed outer-loop quantities, namely the anchor point $\XB_\ell$, the diagonal scaling matrix $\Dx$, and the scalar regularization weight $\cx$, the inner loop approximately minimizes the proximal merit function~\eqref{eq:proximal-merit}. The outer loop then uses the result of this inner minimization to update the anchor point, the diagonal scaling, and the regularization weight.

\subsection{Inner Newton Iteration}
\label{sec:newton}

The first-order optimality condition for the inner problem is
\begin{equation}
	\bm{F}_{\ell}(\X) \triangleq \nabla\mathcal{M}_{\ell}(\X) = \nabla f(\X) + \cx\Dx(\X-\XB_\ell) = \bm{0} .
	\label{eq:inner-residual}
\end{equation}
Newton's method is applied to~\eqref{eq:inner-residual}. At a current inner iterate $\X$, the Newton
correction $\delta\X$ is computed by solving
\begin{equation}
	\bm{J}_{\ell}(\X)\,\delta\X = -\bm{F}_{\ell}(\X),
	\quad \text{where} \quad
	\bm{J}_{\ell}(\X) = \nabla^2\mathcal{M}_{\ell}(\X) = \nabla^2 f(\X) + \cx\Dx .
	\label{eq:inner-newton-system}
\end{equation}
The Newton correction is a descent direction for the merit function $\mathcal{M}_{\ell}(\X)$ whenever $\bm{J}_{\ell}(\X)$ is symmetric positive definite. Indeed, for every nonzero correction $\delta\X$ satisfying~\eqref{eq:inner-newton-system} we have
\begin{equation}
	\nabla\mathcal{M}_{\ell}(\X)\cdot\delta\X = \bm{F}_{\ell}(\X)\cdot\delta\X = -\delta\X\cdot\bm{J}_{\ell}(\X)\delta\X < 0 .
	\label{eq:descent-property}
\end{equation}
Thus, increasing $\cx$, or modifying the diagonal scaling $\Dx$, can improve the conditioning of the inner Hessian and, when necessary, enforce positive definiteness.

\subsubsection{Fallback for a Non-Descent Direction}

If the Newton correction is not a descent direction, this usually indicates that $\bm{J}_{\ell}(\X)$ is not positive definite. In this case, the inner iteration may be terminated with a failure flag. In practice, however, it is often preferable to replace the Newton correction by a scaled gradient direction
\begin{equation}
	\delta\X = -\bm{G}(\X)\bm{F}_{\ell}(\X),
	\label{eq:scaled-gradient-fallback}
\end{equation}
where $\bm{G}(\X)$ is a diagonal scaling matrix defined by
\begin{equation}
	G_{ii}(\X) = \dfrac{1}{
		\max\left\{
		\left|J_{\ell,ii}(\X)\right|,
		\epsilon \left\|\bm{J}_{\ell}(\X)\right\|_{\infty}
		\right\}},
	\quad \text{with} \quad
	\epsilon = \sqrt{\epsilon_{\rm mach}} .
	\label{eq:scaled-gradient-matrix}
\end{equation}
and where $\epsilon_{\rm mach}$ denotes the machine precision (in double precision $\epsilon_{\rm mach}\approx 2.2\cdot 10^{-16}$ and therefore $\epsilon\approx 1.5\cdot 10^{-8}$). Here, the scaling matrix $\bm{G}(\X)$ prevents excessively large components in the fallback direction when the diagonal entries of $\bm{J}_{\ell}(\X)$ are small or when the matrix is nearly singular. This heuristic avoids increasing $\cx$ too aggressively in ill-conditioned cases, while still producing a direction aligned with the negative gradient of the proximal merit function.

\subsubsection{Backtracking Line Search}

To globalize the inner Newton iteration, the full Newton step is combined with a backtracking line search. Given the Newton correction $\delta\X$, the step length $s$ is chosen so that the Armijo condition
\begin{equation}
	\mathcal{M}_{\ell}(\X+s\,\delta\X) \le \mathcal{M}_{\ell}(\X) + c_1\, \nabla\mathcal{M}_{\ell}(\X)\cdot(s\,\delta\X)
	\label{eq:armijo}
\end{equation}
is satisfied. In the reference implementation, the line search starts from $s=1$ and repeatedly replaces $s \gets \beta\, s$, until~\eqref{eq:armijo} holds or a minimum admissible step length is reached. Typical values for the Armijo parameter are in Table~\ref{tab:values}.

\subsubsection{Stagnation, Stalled Iterations, and Small Steps}

During the inner solution, further Newton or line-search iterations may become unproductive. This typically signals that the current local model, scaling or regularization is no longer adequate for the inner problem. We therefore use three diagnostic early-exit tests, none of which is a convergence test.
\begin{itemize}
	\item A \emph{stalled} iteration occurs when the residual is small in a local scale but the current search direction has essentially no first-order effect on the merit function.
	\item \emph{Stagnation} occurs when the residual fails to decrease sufficiently over several consecutive inner iterations.
	\item \emph{Small steps} occur when the accepted corrections are repeatedly negligible relative to the scale of the current iterate. In all cases, when the inner solution process is terminated the control is given back to the outer loop, which may update the proximal weight, modify the diagonal scaling or restart with a more stable local model.
\end{itemize}

Let
\begin{equation*}
	\bm{J}_{\ell}(\X)=\nabla^2\mathcal{M}_{\ell}(\X)
\end{equation*}
be the Hessian of the proximal merit function. We define the diagonal scaling
\begin{equation}
	\bm{S}(\X)=\operatorname{diag}(\bm{w})
	\quad \text{with} \quad
	w_i = \max\left\{
	1,\,
	\sum\nolimits_{j=1}^{n}
	\left|
	\left(\bm{J}_{\ell}(\X)\right)_{ij}
	\right|
	\max\{1,|\X_j|\}
	\right\}.
	\label{eq:stagnation-scaling}
\end{equation}
Thus, $w_i$ estimates the natural local scale of the $i$-th component of the stationarity residual. Using $\bm{J}_{\ell}(\X)$, rather than only $\nabla^2 f(\X)$, makes the scaling consistent with the regularized inner problem.

For an inner iterate $\X_k$, set
\begin{equation}
	r_k = \left\|
	\bm{F}_{\ell}(\X_k)
	\right\|_1,
	\quad \text{and} \quad
	\widehat r_k = \left\|
	\bm{S}(\X_k)^{-1}\bm{F}_{\ell}(\X_k)
	\right\|_1 .
	\label{eq:raw-and-scaled-residuals}
\end{equation}
Here $r_k$ is the unscaled stationarity residual, while $\widehat r_k$ is the
same residual measured in the local scale~\eqref{eq:stagnation-scaling}.

\begin{description}
	\item[Stalled iteration.]
	Let $\delta\X_k$ be the search direction at $\X_k$, and define
	\[
	d_k = \nabla\mathcal{M}_{\ell}(\X_k)\cdot\delta\X_k .
	\]
	The iteration is declared stalled if
	\begin{equation}
		\widehat r_k \leq \tau_{\mathrm{raw}}
		\quad \text{and} \quad
		|d_k| \leq \tau_{\mathrm{slope}} .
		\label{eq:stalled-scaled-residual}
	\end{equation}
	The first condition states that the residual is small in the local curvature-based scale; the second states that the direction has almost zero predicted first-order effect on the merit function. Further backtracking or inner iterations are then unlikely to yield a reliable decrease.
	
	\item[Stagnation.]
	Stagnation is detected from the unscaled residual. A step is productive if
	\begin{equation}
		r_{k+1} \leq \gamma r_k \quad \text{with} \quad 0<\gamma<1 .
		\label{eq:productive-residual-reduction}
	\end{equation}
	If this test fails for $N_{\mathrm{stagnation}}$ consecutive inner iterations while
	\begin{equation}
		\widehat r_k \leq \tau_{\mathrm{scaled}},
		\label{eq:stagnation-scaled-residual}
	\end{equation}
	stagnation is declared: the residual is locally small, but progress toward the inner optimality conditions is insufficient.
	
	\item[Small steps.]
	Let $s_k$ be the line-search step length and $\Delta\X_k=s_k\delta\X_k$ the accepted correction. A small step satisfies
	\begin{equation}
		\|\Delta\X_k\|_\infty \leq \tau_{\mathrm{small}} \max\left\{1,\|\X_k\|_\infty\right\}.
		\label{eq:small-step-test}
	\end{equation}
	If this occurs for $N_{\mathrm{small}}$ consecutive inner iterations, the inner solution process is terminated, since the accepted corrections are negligible at the scale of the current iterate.
\end{description}

\begin{algorithm}[!htb]
	\small
	\caption{Inner iterations}
	\label{alg:algo-inner}
	\begin{algorithmic}[1]
		\State Let $\mathcal{M}_\ell(\X)$ and $\Z_0=\X_\ell$ assigned
		\For{$k = 0,1,2,\ldots,\textrm{\{maximum number of inner iterations\}}$}
		\If{$\|\bm{F}_{\ell}(\Z_k)\|\leq \tau_{\mathrm{inner}}$}
		\State \textsc{Converged}, \Return $\Z_k$
		\EndIf
		\State Compute Newton search direction $\delta\Z$ (Section~\ref{sec:newton}) by solving
		\begin{equation*}
			\bm{J}_{\ell}(\Z_k)\delta\Z
			=
			-\bm{F}_{\ell}(\Z_k),
			\quad \text{where} \quad
			\bm{J}_{\ell}(\Z_k)= \nabla^2 f(\Z_k)+\cx\Dx
		\end{equation*}
		\State Check if $\delta\Z$ is a descent direction, \emph{i.e.}, $\delta\Z^\top(\nabla^2 f(\Z_k)+\cx\Dx)\delta\Z<0$
		
		\If{$\bm{F}_{\ell}(\Z_k)\cdot\delta\Z\leq 0$}
		\State Not a descent direction, use scaled gradient direction \eqref{eq:scaled-gradient-fallback} instead
		\EndIf
		\State $s\leftarrow 1$ \quad $\textrm{found}\leftarrow\textbf{false}$
		\For{$j = 0,1,2,\ldots,\textrm{\{maximum number of backtrackings\}}$}
		\If{ $\mathcal{M}_{\ell}(\Z_k+s\,\delta\Z)\le\mathcal{M}_{\ell}(\Z_k) + c_1\, \nabla\mathcal{M}_{\ell}(\Z_k)\cdot(s\,\delta\Z)$}
		\State Armijo condition satisfied, accept step: $\Z_{k+1}=\Z_k+s\,\delta\Z$
		\State $\textrm{found}\leftarrow\textbf{true}$
		\State \textbf{break}
		\EndIf
		\If{ $s < s_{\min}$ }
		\State \textsc{Failed line-search}. Exit with failure.
		\EndIf
		\State Reduce step: $s\leftarrow \beta s$
		\EndFor
		\If{\textbf{not} found}
		\If{\textsc{stalled}}
		\State \Return last accepted $\Z_k$
		\Else
		\State \textsc{Failed line-search}. Exit with failure.
		\EndIf
		\EndIf
		\State $\Z_{k+1} \leftarrow \Z_k+s\,\delta\Z$
		\If{\textsc{stagnated}}
		\State \Return last accepted $\Z_{k+1}$
		\EndIf
		\EndFor
		\State Reached maximum iteration limits, \Return last accepted $\Z_k$
	\end{algorithmic}
\end{algorithm}

\subsection{Inner Stopping Criteria}

The inner loop terminates as soon as one of the following conditions is met.

\begin{description}
	\item \textbf{Convergence criteria.}
	The inner iteration is successful if at least one stationarity test holds:
	\[
	\frac{1}{n}
	\left\|
	\nabla\mathcal{M}_\ell(\X)
	\right\|_1
	<
	\tau_{\mathrm{inner}},
	\qquad
	\frac{1}{n}
	\left\|
	\bm{S}(\X)\nabla\mathcal{M}_\ell(\X)
	\right\|_1
	<
	\tau_{\mathrm{scaled}} .
	\]
	Here $\bm{S}(\X)$ is a diagonal matrix with
	\[
	S_{ii}(\X)
	=
	\frac{1}{
		\max\left\{
		1,\,
		\left|
		\left[
		\nabla^2\mathcal{M}_\ell(\X)\,\X
		\right]_i
		\right|
		\right\}},
	\]
	so the second test measures stationarity after local curvature-based scaling.
	
	\item \textbf{Interrupt criteria.}
	The inner iteration is unsuccessful but the last computed value is accepted if stagnation, a stalled iteration, small steps or the iteration limit occurs. The latter means that the prescribed maximum $N_{\max}$ of inner iterations is exceeded.
	
	\item \textbf{Failure criteria.}
	The inner iteration is unsuccessful if, before convergence, the line search  fails because $s<s_{\min}$ or the Newton system
	\eqref{eq:inner-newton-system} cannot be solved reliably, for instance due to singularity or severe ill-conditioning.
\end{description}

The returned status is used by the outer loop to decide whether to accept, reject or recompute the trial point with a modified regularization weight.

\subsection{Outer Stopping Criteria and Parameter Updates}
\label{sec:check:convergence}

At the end of each outer iteration, the algorithm checks the prescribed
termination criteria. If none is met, the anchor point and the scalar
regularization weight $\cx$ are updated according to the outcome of the inner Newton solver.

The outer loop terminates successfully when the original objective satisfies
\begin{equation}
	\frac{1}{n}\|\nabla f(\X_{\ell+1})\|_1
	<
	\tau_{\mathrm{outer}},
	\label{eq:outer-stationarity}
\end{equation}
and either the inner solver reports \textbf{successful convergence}, or the \textbf{scaled convergence} criterion holds for $N_{\mathrm{scaled}}$ consecutive outer iterations. This persistence requirement prevents premature termination due to isolated small residuals. The algorithm terminates with a stagnation flag if the inner solver reports \textbf{stagnation} for $N_{\mathrm{stagnation}}$ consecutive outer iterations, since further anchor or regularization updates are then considered unlikely to produce meaningful progress. Otherwise, the iteration continues, and the regularization parameters are updated.

The anchor update follows the accept/reject decision. If the trial point is accepted, $\XB_{\ell+1} \gets \X_{\ell+1}$. If it is rejected, $\XB_{\ell+1} \gets \X_{\ell}$, and the next proximal subproblem is solved with an adjusted $\cx$. Thus, an unsuccessful inner solution does not move the proximal center, while the algorithm can still recover by increasing the regularization strength.

\subsection{Weight Update After an Inner Solution}
\label{sec:cx}

After an inner solution is terminated, the scalar regularization weight $\cx$ is updated multiplicatively according to the progress achieved by the inner iterations. The update has the form
\begin{equation}
	\cx
	=
	\max\{
	\cxmin
	\min\{
	\cxmax,
	\gamma_\ell\, \cx
	\}
	\}
	\label{eq:cx-update}
\end{equation}
that projects the value onto the admissible range for $\cx$, and $\gamma_\ell>0$ is selected from the inner-iteration outcome.

The decision is based on two quantities: the normalized merit variation
\begin{equation}
	\eta_\ell
	\triangleq
	\frac{
		\mathcal{M}_{\ell}(\X_{\ell+1})
		-
		\mathcal{M}_{\ell}(\X_{\ell})
	}{
		\max\left\{
		1,\,
		|\mathcal{M}_{\ell}(\X_{\ell+1})|,\,
		|\mathcal{M}_{\ell}(\X_{\ell})|
		\right\}
	},
	\label{eq:relative-merit-variation}
\end{equation}
and the residual reduction ratio
\begin{equation}
	\rho_\ell
	\triangleq
	\frac{
		\|\nabla\mathcal{M}_{\ell}(\X_{\ell+1})\|_1
	}{
		\|\nabla\mathcal{M}_{\ell}(\X_{\ell})\|_1
	} .
	\label{eq:residual-reduction-ratio}
\end{equation}
Thus, $\eta_\ell<0$ indicates a decrease of the proximal merit function, whereas $\rho_\ell<1$ indicates a reduction of the stationarity residual.

For a fully satisfactory inner solution, $\cx$ is decreased. More precisely, if the inner iteration terminates with the prescribed convergence criterion, then
\begin{equation}
	\gamma_\ell = \theta_{\mathrm{fast}}^{p_\ell},
	\qquad
	p_\ell =
	\frac{n_\ell}{\nu_{\mathrm{rep}}+n_\ell},
	\label{eq:cx-fast-decrease}
\end{equation}
where $n_\ell$ is the number of consecutive occurrences of the same successful termination pattern. Hence, repeated successful inner solutions make the decrease of $\cx$ approach the factor $\theta_{\mathrm{fast}}$. If the inner solution is acceptable but does not satisfy the strongest termination test, the regularization is decreased more conservatively, \emph{i.e.}, by setting
\begin{equation}
	\gamma_\ell = \theta_{\mathrm{slow}} .
	\label{eq:cx-slow-decrease}
\end{equation}
The same conservative decrease is also used when the inner iteration terminates because further progress has become limited, provided that the merit function has decreased and the stationarity residual has been reduced below the prescribed threshold, namely $\eta_\ell < 0$ and $\rho_\ell < \tau_{\rho}$. In this case the regularization is judged effective, although the inner solver has not reached the most stringent stopping condition.

Conversely, if the inner iteration terminates with limited progress and no merit decrease has been obtained, the regularization is increased by setting
\begin{equation}
	\gamma_\ell = \theta_{\mathrm{inc}} .
	\label{eq:cx-slow-increase}
\end{equation}
This strengthens the proximal term in the next subproblem. In all remaining successful or noncritical cases the weight is left unchanged, that is $\gamma_\ell=1$.

In essence, the rule decreases $\cx$ when the inner solution gives reliable evidence that the proximal model is well-behaved, and increases $\cx$ only when the inner iteration stalls without producing a decrease of the merit function. The clipping in \eqref{eq:cx-update} prevents excessively small or excessively large regularization weights.

\subsection{Weight Update After a Failed Inner Solve}
\label{sec:cx:reject}

When the inner iteration is unsuccessful, the scalar regularization weight $\cx$ is updated according to the failure type and to the information collected during the attempted solution. The update distinguishes an inner solution stopped prematurely but still making progress from one failing because the current proximal model is insufficiently regularized.

We distinguish between \emph{interruptions} and \emph{failures}.
\begin{description}
	\item \textbf{Interruption criteria.}
	An interruption stops the inner iteration before convergence, but without a genuine numerical failure; for instance, the maximum number of inner iterations may be reached during the first outer iterations. The indicators $\eta_\ell$ and $\rho_\ell$ are then still used to assess whether the interrupted solution was productive.
	
	If $\eta_\ell<0$, the proximal merit function decreased. If also
	$\rho_\ell<\rho_{\mathrm{big}}$, the stationarity residual was significantly reduced; the solution is therefore considered productive, although not fully converged, and $\cx$ is slowly decreased. The next proximal problem is then less restrictive and may adapt faster to the current anchor point. If instead $\eta_\ell<0$ but $\rho_\ell\geq\rho_{\mathrm{big}}$, the merit decreased without sufficient residual improvement, and $\cx$ is left unchanged.
	
	If $\eta_\ell\geq0$, the inner solver did not decrease the proximal merit function. The regularization is therefore considered too weak for the current local model, and $\cx$ is slowly increased.
	
	\item \textbf{Failure criteria.} A failure corresponds to a genuine numerical difficulty, such as line-search failure, stagnation, an unreliable linear solve, or a non-descent Newton correction that cannot be safely recovered. The default action is to increase $\cx$, thereby strengthening the diagonal regularization, improving the conditioning of the inner Hessian, and stabilizing the next inner problem.
	
	The increase may depend on $\eta_\ell$ and $\rho_\ell$: if $\eta_\ell \geq 0$ or $\rho_\ell \geq \rho_{\mathrm{big}}$, $\cx$ is increased more aggressively. If the failure occurs after a partial merit decrease and a moderate residual reduction, a milder increase is used.
\end{description}

\section{Numerical Tests}

Our experimental setup consists of an implementation of Algorithms~\ref{alg:algo-base} and~\ref{alg:algo-inner} with the parameter values in Table~\ref{tab:values} and the weight update mechanism described above. We tested the proposed scheme on 100 benchmark problems from the literature.

\begin{table}[!htb]
	\centering
	\caption{Parameter configuration used in the numerical experiments.}
	\label{tab:values}
	\begin{minipage}[t]{0.47\textwidth}
		\centering
		\begin{tabular}{lc}
			\toprule
			\multicolumn{2}{c}{Tolerances and iteration limits} \\
			\midrule
			$\tau_{\mathrm{inner}}$      & $10^{-12}$ \\
			$\tau_{\mathrm{outer}}$      & $10^{-10}$ \\
			$\tau_{\mathrm{scaled}}$     & $10^{-16}$ \\
			$\tau_{\mathrm{stagnation}}$ & $10^{-12}$ \\
			$\tau_{\rho}$                & $10^{-1}$ \\
			$\tau_{\mathrm{raw}}$        & $10^{-6}$ \\
			\addlinespace
			$N_{\max}$                   & $20$ \\
			$N_{\mathrm{scaled}}$        & $5$ \\
			$N_{\mathrm{stagnation}}$    & $3$ \\
			\bottomrule
		\end{tabular}
	\end{minipage}
	\hfill
	\begin{minipage}[t]{0.47\textwidth}
		\centering
		\begin{tabular}{lc}
			\toprule
			\multicolumn{2}{c}{Algorithmic parameters} \\
			\midrule
			$s_{\min}$                 & $10^{-20}$ \\
			$\gamma$                   & $0.5$ \\
			$\beta$                    & $0.75$ \\
			$c_1$                      & $10^{-4}$ \\
			$\nu_{\mathrm{rep}}$       & $0.25$ \\
			$\theta_{\mathrm{fast}}$   & $0.5$ \\
			$\theta_{\mathrm{slow}}$   & $0.8$ \\
			$\theta_{\mathrm{inc}}$    & $1.5$ \\
			$\cx_{\min}$               & $10^{-18}$ \\
			$\cx_{\max}$               & $10^{10}$ \\
			\bottomrule
		\end{tabular}
	\end{minipage}
\end{table}

\subsection{Convergence Criteria}

For all tests, the initial regularization weight was set to $10^{-18}$,
which applies a mild initial proximal penalty to allow larger initial
Newton steps if the landscape permits. A run is considered successfully converged if the solver successfully terminates and the infinity norm of the objective function gradient at the final iterate satisfies the strict absolute tolerance condition:
\[
\|\nabla f(x)\|_\infty \le \num[round-precision = 0]{1e-12}.
\]
or a scaled relative one
\[
\dfrac{\|\nabla f(x)\|_\infty}{\max\left\{1,\|\nabla^2 f(x)\|_\infty\right\}} \le  \num[round-precision = 0]{1e-12}.
\]
Runs where the solver stops (\emph{e.g.}, due to step stagnation) but fails to meet $\num[round-precision = 0]{1e-4}$ absolute tolerance are classified as ``fake converged''.

\subsection{Results Summary}

The algorithm demonstrated good robustness across the benchmark suite. Out of the 100 problems tested, the solver reported convergence for all of them, with zero complete failures. Applying the strict gradient tolerance, 98 problems achieved true convergence. Table~\ref{tab:global_summary} presents the aggregated statistics for the benchmark. For the 98 fully successful runs, the algorithm required an average of $16.4$ outer iterations and $228.45$ inner Newton steps. This indicates that the outer loop effectively adapts the anchor point and the diagonal scaling, keeping the required sub-problem iterations manageable.

\begin{table}[!htb]
	\centering
	\caption{Performance summary of the proposed solver over 100 benchmark problems.}
	\label{tab:global_summary}
	\begin{tabular}{@{}lc@{}}
		\midrule
		Total Benchmark Problems & 100 \\
		Solver-Reported Converged & 100 \\
		\textbf{Accepted Converged} ($\|\nabla f\|_\infty \le 10^{-4}$) & \textbf{98} \\
		Fake Converged & 2 \\
		Failed Runs & 0 \\
		\midrule
		Average Outer Iterations (Accepted) & 16.40 \\
		Average Sub-iterations (Accepted) & 228.45 \\
		Total Outer Iterations (Accepted) & 1607 \\
		Total Sub-iterations (Accepted) & 22388 \\
		\midrule
	\end{tabular}
\end{table}

\subsection{Analysis of Challenging Problems}

Only two problems in the suite, namely the Extended Helical Valley ($n=100$) and the Meyer growth-model fit ($n=3$), triggered a fake convergence status. The details of these specific runs are reported in Table~\ref{tab:fake_converged}.

\begin{table}[htbp]
	\centering
	\caption{Details of the ``fake converged'' problems that stagnated just above the strict gradient tolerance.}
	\label{tab:fake_converged}
	\begin{tabular}{@{~}l@{~~}c@{~~}c@{~~}c@{~~}c@{~~}c@{~}}
		\toprule
		No. & $n$ & Outer & Inner & Final $f(x)$ & Final $\|\nabla f(\X)\|_\infty$ \\
		\midrule
		52 & 100 & 45 & 832 & $\num{9.83e-9}$ & $\num{1.98e-4}$ \\
		94 & 3   & 62 & 1001 & $87.94$ & $\num{9.83e-3}$ \\
		\bottomrule
	\end{tabular}
\end{table}

In both cases, the algorithm drastically reduced the objective function and reached the immediate vicinity of the minimum but ultimately stagnated. For test 52, the gradient norm ($1.98 \times 10^{-4}$) missed the target threshold by a negligible margin. The Meyer function (test 94) is notoriously highly nonlinear and poorly scaled, which routinely causes premature termination in classical unconstrained solvers due to flat valleys. The fact that the dynamic proximal point method avoids divergence and safely halts in these extreme cases confirms the safety and stability of the adaptive fallbacks implemented in the inner loop.

\begingroup
\scriptsize
\setlength{\LTleft}{\fill}
\setlength{\LTright}{\fill}
\begin{longtable}{@{~}l@{~}p{4cm}@{~}c@{~}c@{~}c@{~}c@{~}c@{~}c@{~}}
	\caption{Statistics with benchmark references.}
	\label{tab:proxi-minimize-low-with-refs}\\
	\toprule
	No. & Test & $n$ & Outer & Inner & Final $f(\X)$ & Final $\|\nabla f(\X)\|_\infty$ & Ref. \\
	\midrule
	\endfirsthead
	\toprule
	No. & Test & $n$ & Iter & SubIter & Final $f(\X)$ & Final $\|\nabla f(\X)\|_\infty$ & Ref. \\
	\midrule
	\endhead
	\midrule
	\multicolumn{8}{r}{Continued on next page}\\
	\endfoot
	\bottomrule
	\endlastfoot
	1  & Arrowhead quartic                & 100 & 25  & 25   & \num{0.0e0}         & \num{2.665e-14} & \cite{Andrei2008} \\
	2  & Arrowhead quartic                & 200 & 25  & 25   & \num{0.0e0}         & \num{2.665e-14} & \cite{Andrei2008} \\
	3  & Brown almost-linear              & 100 & 4   & 21   & \num{1.26218e-29}   & \num{7.105e-15} & \cite{Brown1969} \\
	4  & Brown almost-linear              & 200 & 7   & 56   & \num{1.63174e-25}   & \num{1.148e-11} & \cite{Brown1969} \\
	5  & Discrete integral equation       & 100 & 3   & 5    & \num{4.03675e-30}   & \num{1.036e-15} & \cite{More1981TUO} \\
	6  & Discrete integral equation       & 200 & 3   & 5    & \num{7.11581e-30}   & \num{9.766e-16} & \cite{More1981TUO} \\
	7  & Dixon-Price chain                & 100 & 3   & 9    & \num{9.18283e-30}   & \num{3.553e-14} & \cite{DixonPrice1989} \\
	8  & Dixon-Price chain                & 200 & 3   & 9    & \num{9.18283e-30}   & \num{3.553e-14} & \cite{DixonPrice1989} \\
	9  & Diagonal 1 function              & 100 & 3   & 9    & \num{-1.5706e4}     & \num{4.263e-14} & \cite{Andrei2008} \\
	10 & Diagonal 1 function              & 200 & 3   & 10   & \num{-7.6396e4}     & \num{8.527e-14} & \cite{Andrei2008} \\
	11 & Diagonal 2 function              & 100 & 3   & 11   & \num{1.57414e1}     & \num{7.147e-16} & \cite{Andrei2008} \\
	12 & Diagonal 2 function              & 200 & 3   & 12   & \num{1.98545e1}     & \num{2.776e-17} & \cite{Andrei2008} \\
	13 & Diagonal 3 function              & 100 & 3   & 7    & \num{-4.6058e3}     & \num{9.770e-15} & \cite{Andrei2008} \\
	14 & Diagonal 3 function              & 200 & 3   & 7    & \num{-1.91825e4}    & \num{1.954e-14} & \cite{Andrei2008} \\
	15 & Diagonal 4 quadratic             & 100 & 3   & 2    & \num{0e0}           & \num{0.000e+00} & \cite{Andrei2008} \\
	16 & Diagonal 4 quadratic             & 200 & 3   & 2    & \num{0e0}           & \num{0.000e+00} & \cite{Andrei2008} \\
	17 & Diagonal 5 function              & 100 & 3   & 6    & \num{6.93147e1}     & \num{2.346e-13} & \cite{Andrei2008} \\
	18 & Diagonal 5 function              & 200 & 3   & 6    & \num{1.38629e2}     & \num{2.346e-13} & \cite{Andrei2008} \\
	19 & Extended Beale blocks            & 82  & 3   & 11   & \num{1.81931e-29}   & \num{5.995e-15} & \cite{Andrei2008} \\
	20 & Extended Beale blocks            & 200 & 3   & 11   & \num{4.43734e-29}   & \num{5.995e-15} & \cite{Andrei2008} \\
	21 & Extended Hiebert                 & 100 & 360 & 6421 & \num{4.77261e-25}   & \num{1.954e-13} & \cite{Andrei2008} \\
	22 & Extended Hiebert                 & 200 & 133 & 2119 & \num{0e0}           & \num{0.000e+00} & \cite{Andrei2008} \\
	23 & Extended quadratic penalty 1     & 100 & 5   & 18   & \num{3.90062e2}     & \num{9.275e-12} & \cite{Andrei2008} \\
	24 & Extended quadratic penalty 1     & 200 & 3   & 11   & \num{7.90031e2}     & \num{5.289e-12} & \cite{Andrei2008} \\
	25 & Extended quadratic penalty 2     & 100 & 7   & 69   & \num{5.40402e-2}    & \num{4.282e-13} & \cite{Andrei2008} \\
	26 & Extended quadratic penalty 2     & 200 & 7   & 72   & \num{2.16161e-1}    & \num{1.611e-12} & \cite{Andrei2008} \\
	27 & Expanded Schaffer F6             & 100 & 3   & 9    & \num{9.71591e-1}    & \num{1.541e-13} & \cite{Schwefel1995} \\
	28 & Expanded Schaffer F6             & 200 & 3   & 9    & \num{1.94318e0}     & \num{1.535e-13} & \cite{Schwefel1995} \\
	29 & Extended BD1 block diagonal      & 100 & 6   & 68   & \num{0e0}           & \num{0.000e+00} & \cite{Andrei2008} \\
	30 & Extended BD1 block diagonal      & 200 & 6   & 68   & \num{0e0}           & \num{0.000e+00} & \cite{Andrei2008} \\
	31 & Extended Freudenstein-Roth       & 100 & 5   & 7    & \num{2.44921e3}     & \num{3.771e-07} & \cite{Andrei2008} \\
	32 & Extended Freudenstein-Roth       & 200 & 5   & 19   & \num{4.89843e3}     & \num{4.236e-09} & \cite{Andrei2008} \\
	33 & Extended TET function            & 100 & 23  & 219  & \num{1.27963e2}     & \num{0.000e+00} & \cite{Andrei2008} \\
	34 & Extended TET function            & 200 & 23  & 219  & \num{2.55927e2}     & \num{0.000e+00} & \cite{Andrei2008} \\
	35 & Extended Wood function           & 100 & 3   & 11   & \num{1.9695e2}      & \num{3.286e-14} & \cite{Colville1968} \\
	36 & Extended Wood function           & 200 & 3   & 11   & \num{3.93899e2}     & \num{3.286e-14} & \cite{Colville1968} \\
	37 & Full Hessian FH1                 & 100 & 3   & 3    & \num{8.92008e-28}   & \num{1.098e-13} & \cite{Andrei2008} \\
	38 & Full Hessian FH1                 & 200 & 3   & 3    & \num{1.41675e-26}   & \num{2.444e-13} & \cite{Andrei2008} \\
	39 & Full Hessian FH2                 & 100 & 3   & 2    & \num{7.24642e-23}   & \num{4.545e-12} & \cite{Andrei2008} \\
	40 & Full Hessian FH2                 & 200 & 3   & 3    & \num{3.78191e-26}   & \num{4.405e-13} & \cite{Andrei2008} \\
	41 & Generalized tridiagonal 1        & 100 & 3   & 17   & \num{9.72103e1}     & \num{4.441e-15} & \cite{Andrei2008} \\
	42 & Generalized tridiagonal 1        & 200 & 3   & 19   & \num{1.9721e2}      & \num{1.090e-13} & \cite{Andrei2008} \\
	43 & Griewank product-sum             & 100 & 3   & 6    & \num{0e0}           & \num{1.000e-18} & \cite{Griewank1981} \\
	44 & Griewank product-sum             & 200 & 3   & 6    & \num{0e0}           & \num{1.000e-18} & \cite{Griewank1981} \\
	45 & Generalized Rosenbrock           & 100 & 47  & 812  & \num{4.19365e-24}   & \num{2.376e-12} & \cite{Andrei2008} \\
	46 & Generalized Rosenbrock           & 200 & 60  & 1020 & \num{5.12644e-23}   & \num{8.238e-12} & \cite{Andrei2008} \\
	47 & Generalized White-Holst chain    & 100 & 72  & 1260 & \num{1.16225e-20}   & \num{7.616e-11} & \cite{Andrei2008} \\
	48 & Generalized White-Holst chain    & 200 & 104 & 2002 & \num{1.8588e-24}    & \num{9.992e-13} & \cite{Andrei2008} \\
	49 & Hager function                   & 100 & 3   & 7    & \num{-6.53079e2}    & \num{3.730e-14} & \cite{Andrei2008} \\
	50 & Hager function                   & 200 & 3   & 8    & \num{-2.4932e3}     & \num{3.553e-15} & \cite{Andrei2008} \\
	51 & Fletcher-Powell Helical Valley   & 3   & 3   & 21   & \num{4.29415e-42}   & \num{1.002e-20} & \cite{FletcherPowell1963} \\
	52 & Extended Fletcher-Powell         & 100 & 52  & 832  & \num{9.83924e-09}   & \color{red}\num{1.984e-04} & \cite{FletcherPowell1963}\\
	53 & Levy function                    & 100 & 41  & 662  & \num{1.01869e2}     & \num{1.315e-11} & \cite{LevyMontalvo1985} \\
	54 & Levy function                    & 200 & 34  & 568  & \num{2.05777e2}     & \num{2.708e-11} & \cite{LevyMontalvo1985} \\
	55 & Extended Powell singular func.   & 100 & 4   & 29   & \num{3.86299e-16}   & \num{1.745e-12} & \cite{Powell1962} \\
	56 & Extended Powell singular func.   & 200 & 4   & 29   & \num{7.72598e-16}   & \num{1.745e-12} & \cite{Powell1962} \\
	57 & Perturbed quadratic              & 100 & 3   & 5    & \num{1.49326e-52}   & \num{6.006e-26} & \cite{Andrei2008} \\
	58 & Perturbed quadratic              & 200 & 3   & 5    & \num{8.49437e-52}   & \num{1.334e-25} & \cite{Andrei2008} \\
	59 & Perturbed quadratic diagonal     & 100 & 3   & 2    & \num{4.68828e-26}   & \num{6.366e-14} & \cite{Andrei2008} \\
	60 & Perturbed quadratic diagonal     & 200 & 3   & 2    & \num{1.51917e-25}   & \num{1.510e-13} & \cite{Andrei2008} \\
	61 & Qing function                    & 100 & 23  & 258  & \num{4.77202e-23}   & \num{2.579e-11} & \cite{HedarTestGO} \\
	62 & Qing function                    & 200 & 54  & 733  & \num{8.15773e-22}   & \num{1.142e-10} & \cite{HedarTestGO} \\
	63 & Quadratic QF1                    & 100 & 3   & 2    & \num{-5e-3}         & \num{8.882e-16} & \cite{Andrei2008} \\
	64 & Quadratic QF1                    & 200 & 3   & 2    & \num{-2.5e-3}       & \num{8.882e-16} & \cite{Andrei2008} \\
	65 & Quadratic QF2                    & 100 & 3   & 7    & \num{-1.00125e0}    & \num{3.517e-13} & \cite{Andrei2008} \\
	66 & Quadratic QF2                    & 200 & 3   & 7    & \num{-1.00062e0}    & \num{1.421e-14} & \cite{Andrei2008} \\
	67 & Rastrigin oscillatory            & 100 & 3   & 13   & \num{1.59192e3}     & \num{4.530e-14} & \cite{Rastrigin1974} \\
	68 & Rastrigin oscillatory            & 200 & 3   & 8    & \num{3.18385e3}     & \num{4.441e-14} & \cite{Rastrigin1974}  \\
	69 & Raydan 1 function                & 100 & 5   & 48   & 505                 & \num{1.066e-15} & \cite{Andrei2008} \\
	70 & Raydan 1 function                & 200 & 5   & 48   & 2010                & \num{2.842e-15} & \cite{Andrei2008} \\
	71 & Raydan 2 function                & 100 & 3   & 7    & 100                 & \num{0.000e+00} & \cite{Andrei2008} \\
	72 & Raydan 2 function                & 200 & 3   & 7    & 200                 & \num{0.000e+00} & \cite{Andrei2008} \\
	73 & Extended Rosenbrock function     & 100 & 42  & 757  & \num{4.47115e-27}   & \num{6.661e-14} & \cite{Rosenbrock1960} \\
	74 & Extended Rosenbrock function     & 200 & 59  & 967  & \num{1.29686e-22}   & \num{1.308e-11} & \cite{Rosenbrock1960} \\
	75 & Schwefel cumulative quadratic    & 100 & 3   & 3    & \num{3.1855e-31}    & \num{9.702e-15} & \cite{Liang2013CEC2014} \\
	76 & Schwefel cumulative quadratic    & 200 & 3   & 3    & \num{2.6071e-30}    & \num{3.940e-14} & \cite{Liang2013CEC2014} \\
	77 & Trid chain                       & 100 & 3   & 2    & \num{-1.716e5}      & \num{4.547e-13} & \cite{Ali2005} \\
	78 & Trid chain                       & 200 & 3   & 2    & \num{-1.3532e+06}   & \num{3.638e-12} &  \cite{Ali2005}\\
	79 & Zakharov function                & 100 & 4   & 27   & \num{3.59957e-28}   & \num{7.462e-13} & \cite{Zhigljavsky1991} \\
	80 & Zakharov function                & 200 & 4   & 32   & \num{4.60736e-31}   & \num{8.675e-16} & \cite{Zhigljavsky1991} \\
	81 & Bard rational fitting            & 3   & 3   & 15   & \num{8.21488e-3}    & \num{8.327e-16} & \cite{More1981TUO} \\
	82 & Biggs EXP6                       & 6   & 44  & 618  & \num{5.65565e-3}    & \num{2.608e-12} & \cite{More1981TUO} \\
	83 & Box 3D exponential fit           & 3   & 3   & 11   & \num{6.77927e-32}   & \num{5.007e-16} & \cite{More1981TUO} \\
	84 & Brown badly scaled LS            & 2   & 3   & 12   & \num{0e0}           & \num{0.000e+00} & \cite{More1981TUO}  \\
	85 & Brown-Dennis quartic LS          & 4   & 5   & 30   & \num{8.58222e4}     & \num{1.310e-10} & \cite{More1981TUO} \\
	86 & Broyden banded LS                & 20  & 3   & 13   & \num{1.88924e1}     & \num{4.774e-15} & \cite{More1981TUO} \\
	87 & Broyden tridiagonal LS           & 20  & 3   & 21   & \num{1.60809}       & \num{1.233e-15} & \cite{More1981TUO} \\
	88 & Discrete boundary-value          & 20  & 3   & 4    & \num{2.36521e-25}   & \num{1.391e-14} & \cite{More1981TUO} \\
	89 & Extended Powell singular         & 12  & 4   & 29   & \num{4.63559e-17}   & \num{1.745e-12} & \cite{More1981TUO} \\
	90 & Extended Rosenbrock valley       & 10  & 32  & 525  & \num{3.33748e-25}   & \num{2.665e-13} & \cite{More1981TUO} \\
	91 & Gaussian fitting                 & 3   & 3   & 4    & \num{1.12793e-08}   & \num{2.466e-16} & \cite{More1981TUO} \\
	92 & Jennrich-Sampson exp. fit        & 2   & 5   & 30   & \num{124.362} & \num{1.137e-12} & \cite{More1981TUO} \\
	93 & Kowalik-Osborne rational fit     & 4   & 33  & 519  & \num{0.000307486} & \num{6.447e-13} & \cite{More1981TUO} \\
	94 & Meyer growth-model fit           & 3   & 62  & 1001 & \num{87.9459} & \color{red}\num{9.831e-03} & \cite{More1981TUO} \\
	95 & Osborne 1 biexponential fit      & 5   & 4   & 29   & \num{5.46489e-05} & \num{6.504e-15} & \cite{More1981TUO} \\
	96 & Penalty I ill-conditioned model  & 10  & 4   & 41   & \num{7.08765e-05} & \num{9.341e-17} & \cite{More1981TUO} \\
	97 & Powell badly scaled LS           & 2   & 63  & 823  & \num{6.66404e-26} & \num{5.163e-13} & \cite{More1981TUO} \\
	98 & Trigonometric system LS          & 10  & 39  & 609  & \num{0.000620651} & \num{9.877e-13} & \cite{More1981TUO} \\
	99 & Variably dimensioned valley      & 10  & 3   & 18   & \num{0} & \num{0.000e+00} & \cite{More1981TUO} \\
	100 & Watson least-square             & 12  & 3   & 14   & \num{4.72238e-10} & \num{1.708e-13} & \cite{More1981TUO} \\
\end{longtable}
\endgroup

\section{Conclusion}

We have presented a dynamic proximal point method for unconstrained
optimization. The method combines an adaptive diagonal regularization with a Newton-type inner solver. The proximal term stabilizes the local model, while the diagonal scaling and the scalar weight are updated dynamically according to the behavior of the inner iteration. The resulting algorithm adjusts the regularization strength automatically. When the inner solver produces a reliable decrease of the merit function and a significant reduction of the stationarity residual, the proximal weight can be relaxed, allowing the method to approach a standard Newton iteration. Conversely, when the inner solver stagnates, fails or encounters an indefinite or poorly conditioned local model, the regularization is strengthened. This mechanism improves robustness on ill-conditioned and non-convex problems, while preserving fast local progress when the Newton model is reliable. Such a formulation is accompanied by explicit stopping criteria, failure handling rules, update strategies, and implementation-oriented parameter choices. All these features make the method suitable for direct implementation and for integration into existing nonlinear optimization codes.

The same framework can be extended to constrained optimization problems. Equality constraints can be incorporated by introducing suitable residual terms, while inequality and box constraints can be handled through slack variables, active-set mechanisms or bound-aware proximal terms. In this setting, the proximal regularization can be used to construct a merit function whose local model can be convexified by increasing the regularization weight. The constrained case is more involved, because feasibility, complementarity, constraint qualification,and multiplier updates must be treated consistently. These extensions will be addressed in future work.

\subsection*{Acknowledgments}
\TheFunding

\bibliographystyle{habbrv}
\bibliography{bibliography}
	
\end{document}